\documentclass[letterpaper, 10 pt, conference]{ieeeconf}
\IEEEoverridecommandlockouts
\usepackage{amssymb}
\usepackage{amsmath}
\usepackage[T1]{fontenc}
\usepackage[latin1]{inputenc}
\usepackage{url}
\usepackage{cite}
\usepackage{amsfonts}
\usepackage{graphicx}
\usepackage{epsfig}

\newtheorem{example1}{Example}[section]

\newtheorem{lemma*}{Lemma}

\begin{document}
\title{\textbf{Bilateral Boundary Control of One-Dimensional First- and Second-Order PDEs using Infinite-Dimensional Backstepping}}

\author{{Rafael Vazquez and Miroslav Krstic}
\thanks{R. Vazquez is with the Department of Aerospace Engineering, Universidad de Sevilla, Camino de los Descubrimiento s.n., 41092 Sevilla, Spain.}%
\thanks{M. Krstic is with the Department of Mechanical and Aerospace Engineering, University of California San Diego, La Jolla, CA 92093-0411, USA.}%
}

\maketitle
 \interdisplaylinepenalty=2500



\baselineskip=.95 \normalbaselineskip

\begin{abstract}
This paper develops an extension of infinite-dimensional backstepping method for parabolic and hyperbolic systems in one spatial dimension with two actuators. Typically, PDE backstepping is applied in 1-D domains with an actuator at one end. Here, we consider the use of two actuators, one at each end of the domain, which we refer to as bilateral control (as opposed to unilateral control). Bilateral control laws are derived for linear reaction-diffusion, wave and $2\times2$ hyperbolic 1-D systems (with same speed of transport in both directions). The extension is nontrivial but straightforward if the backstepping transformation is adequately posed. The resulting bilateral controllers are compared with their unilateral counterparts in the reaction-diffusion case for constant coefficients, by making use of explicit solutions, showing a  reduction in control effort as a tradeoff for the presence of two actuators when the system coefficients are large. These results open the door for more sophisticated designs such as bilateral sensor/actuator output feedback and fault-tolerant designs.
\end{abstract}

\section{Introduction}
The backstepping method has proved itself to be an ubiquitous method for PDE control. First developed to design feedback control laws and observers for one-dimensional reaction-diffusion PDEs~\cite{krstic}, it has since been applied to many other systems including, among others, flow control~\cite{vazquez,vazquez-coron}, thermal loops~\cite{convloop}, nonlinear PDEs~\cite{vazquez2}, hyperbolic 1-D systems~\cite{vazquez-nonlinear,florent,krstic3},  multi-agent deployment~\cite{jie}, wave equations~\cite{krstic2}, and delays~\cite{krstic5}. Some of the more striking features of backstepping include the possibility of finding explicit control laws in some cases (see e.g.~\cite{Vazquez2014}) or even designing adaptive controllers~\cite{krstic4}.

Backstepping has been typically applied to PDEs formulated in 1-D domains, with actuation at one end of the domain. This paper develops backstepping for the case when two actuators are available, one at each end of the domain. We refer to this situation as bilateral control---as opposed to the single actuator case, denoted as unilateral control. The bilateral case is a nontrivial extension of the method that needs a specific formulation of the backstepping transformation to be able to solve the stabilization problem. More concretely, whereas the backstepping transformation in the unilateral case is formulated as an integral starting at the non-actuated end, in the bilateral case the integral starts at the middle of the domain. Based in this simple idea it is straightforward to extend the unilateral results and even obtain explicit results for the constant coefficient case. Explicit solutions allow to easily compare the resulting bilateral controllers. We compare unilateral and bilateral controllers for  the reaction-diffusion equation, and show that, when the system coefficients are large, there is a considerable reduction in control effort as a tradeoff for the presence of two actuators. 

The focus of this paper is on design, stabilization, and finding explicit controllers when available. Closed-loop well-posedness is assumed, in the most basic Sobolev space appropriate for each type of PDE. Due to the linearity of the equations and the good properties of backstepping a detailed well-posedness analysis would be rather straightforward but lengthy and dependent on the type of the equation. We only study the well-posedness of the kernel equations which is critical towards finding the control laws.

The results of this paper were originally inspired by~\cite{nball}, where the problem of stabilization of reaction-diffusion equations in balls of arbitrary dimension is addressed. Since a ball in dimension 1 is actually an interval, and the boundaries are the two ends; then, applying the formulas of~\cite{nball} one obtains the results that will be shown in Section~\ref{sect-explicit}. In this paper the authors explore similar results for other types of equations.

This paper presents bilateral control laws for linear reaction-diffusion, wave and $2\times2$ hyperbolic 1-D systems. For $2\times2$ hyperbolic 1-D systems, only the case of both states having the same speed of transport is considered. The reason to only consider this particular case is that the analysis of the resulting kernel equations is straightforward. It is also the case used in the paper to build the control design for the wave equation. Finally, it allows to derive explicit solutions in the constant-coefficient case. The case of different speeds of transport can also be addressed, however while writing this paper the authors learnt of a new paper~\cite{Fl-bilateral} that already solves this more general and challenging problem.

The paper is organized as follows. In Section~\ref{sect-parabolic} we solve the bilateral control problem for a reaction-diffusion equation. In Section~\ref{sect-hyp} we continue with a bilateral design for $2\times2$ hyperbolic 1-D systems with the same speed of transport. This design is then adapted for wave equations in Section~\ref{sect-wave}. Explicit controllers for all cases are presented in Section~\ref{sect-explicit}. We then compare bilateral and unilateral results in Section~\ref{sect-comparison}. We finish in Section~\ref{sect:conclusions} with some concluding remarks. 

\section{Reaction-Diffusion PDEs}~\label{sect-parabolic}
Consider the reaction-diffusion equation
\begin{eqnarray} \label{eqn-urd}
u_t&=&\epsilon u_{xx}+\lambda(x) u,
\end{eqnarray}
for $t>0$, in the domain\footnote{In PDE backstepping the domain is typically written as $x\in[0,L]$. However for the bilateral formulation it is somewhat convenient to fix the origin at the middle point of the domain.} $x\in[-L/2,L/2]$, with $\epsilon>0$, $\lambda(x)$ a differentiable function, and 
with boundary conditions 
\begin{eqnarray}\label{eqn-urdbc1}
u(t,L)&=&U_1(t),\\
u(t,-L)&=&U_2(t),\label{eqn-urdbc2}
\end{eqnarray}
where $U_1$ and $U_2$ are actuator variables.
For sufficiently large $\lambda(x)>0$, (\ref{eqn-urd})--(\ref{eqn-urdbc2}) is open-loop unstable.

To design feedback control laws for $U_1$ and $U_2$, consider a transformation defined as
\begin{equation}\label{eqn-transf-unidim}
w(t,x)=u(t,x)-\int_{-x}^x K(x,\xi) u(t,\xi)d\xi,
\end{equation}
with $w(t,x)$ (the \emph{target} variable) verifying the following system (\emph{target} system):
\begin{eqnarray}\label{eqn-wpde}
w_t&=&\epsilon w_{xx},\\
w(t,L)&=&w(t,-L)=0.\label{eqn-wbc}
\end{eqnarray}

\begin{figure}[ht]
\includegraphics[width=7cm]{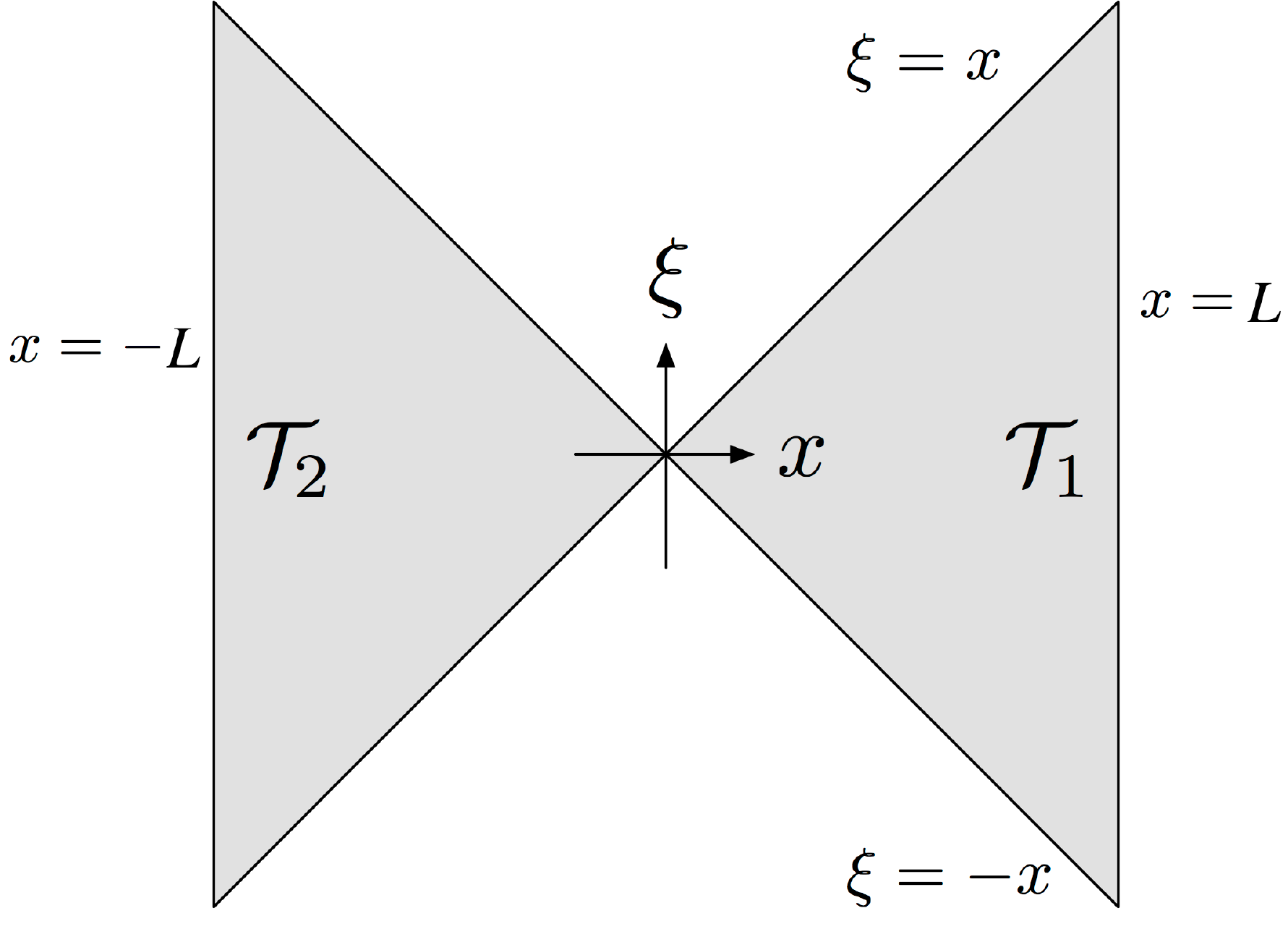}
\centering
\caption{The domain $\mathcal T$ and subdomains $\mathcal T_1$, $\mathcal T_2$.}
\label{fig:hourglass}
\end{figure}

Working out the kernel equations as in the unilateral case (see e.g.~\cite{krstic}), one finds
\begin{eqnarray}
\epsilon K_{xx}(x,\xi)-\epsilon K_{\xi \xi}(x,\xi)&=&
 \lambda(\xi) K(x,\xi) \label{eqn-K1}
,\\
K(x,x)&=&-\int_0^x \frac{\lambda(\xi)}{2\epsilon} d\xi,\\
K(x,-x)&=&0 \label{eqn-K3}
\end{eqnarray}
in the hourglass-shaped domain $\mathcal T=\{(x,\xi):x\in[-L,L],\,\xi \in[-\vert x \vert,\vert x \vert\}$, represented in Fig.~\ref{fig:hourglass}. It is possible to separate the domain into two: $\mathcal T=\mathcal T_1\cup\mathcal T_2$, where $\mathcal T_1=\{(x,\xi):x\in[0,L],-x \leq \xi \leq x\}$ and $\mathcal T_2=\{(x,\xi):x\in[-L,0],x \leq \xi \leq -x\}$, as shown in Fig.~\ref{fig:hourglass}. It is clear that the interior of these domains is disjoint and it follows that the kernel equations can be solved separately in each of the domains. In addition, it is easy to see that if the equations are well-posed in the domain $\mathcal T_1$ they must be in $\mathcal T_2$ by a symmetry argument (switching variables from $(x,\xi)$ to $(\hat x,\hat \xi)=(-x,\xi)$ maps the domain $\mathcal T_2$ into $\mathcal T_1$ and leaves the structure of the equations unchanged, except for some sign switch).
Once the kernel equations are solved, the control laws are given by the values of the kernel at both ends of the domain $\mathcal T$ as follows
\begin{eqnarray}\label{eqn-U1}
U_1&=&\int_{-L}^L K(L,\xi) u(\xi)d\xi,\\
U_2&=&-\int_{-L}^L K(-L,\xi) u(\xi)d\xi.\label{eqn-U2}
\end{eqnarray}
Then, and assuming the backstepping transformation is invertible, the state $u(t,x)$ ``inherits'' the stability properties of the target state $w(t,x)$, whose origin, in view of its defining PDE (\ref{eqn-wpde})--(\ref{eqn-wbc}), is clearly exponentially stable. Now, since inverting the backstepping transformation amounts to solving a Volterra integral equation of the second kind, then inversion is always possible under very mild conditions~\cite{invertibility} (for instance if the kernel is at least bounded, which it is in our case).

Thus, the problem of designing an stabilizing control law is reduced to solving (\ref{eqn-K1})--(\ref{eqn-K3}) so that (\ref{eqn-U1})--(\ref{eqn-U2}) can be implemented and showing it is at least bounded. We first address the question of existence and uniqueness in the domain $\mathcal T_1$. Consider the new variables $\alpha=x+\xi$ and $\beta=x-\xi$ and the kernel $K=G(\alpha,\beta)$ a function of the new variables. Written in terms of $(\alpha,\beta)$, (\ref{eqn-K1})--(\ref{eqn-K3}) become
\begin{eqnarray}\label{eqn-KG1}
G_{\alpha\beta}(\alpha,\beta)&=&\frac{\lambda\left(\frac{\alpha-\beta}{2}\right)}{4\epsilon} G(\alpha,\beta),\\
G(\alpha,0)&=&-\int_0^\frac{\alpha}{2} \frac{\lambda(\xi)}{2\epsilon} d\xi,\\
G(0,\beta)&=&0 \label{eqn-KG3}
\end{eqnarray}
in the domain $\mathcal T_1$ which is now expressed as $\mathcal T_1=\left\{(\alpha,\beta):\alpha>0,\beta>0,\alpha+\beta\leq2L\right\}$.

The problem in this form is known as Goursat's problem~\cite{day}, which interestingly appears in the theory of nonlinear wave propagation~\cite{jeffrey}. It is known that if $\lambda$ is at least differentiable then there is a unique $C^1(\mathcal T_1)$ solution~\cite{holten}. This can also be shown directly by transforming (\ref{eqn-KG1})-(\ref{eqn-KG3}) into an integral equation and then applying the method of successive approximations, as in~\cite{krstic}.

\subsection{An equivalent problem}
The problem of stabilization of systems of reaction-diffusion equations with the same diffusion coefficients was addressed in~\cite{orlov}. Now we show the relationship between that result and the bilateral control problem for reaction-diffusion equations. The first step is to ``fold'' the state $u$ into a $2\times2$ system $(u_1,u_2)$, defining 
\begin{equation}
u(t,x)=\left\{ \begin{array}{ll} u_1(t,x), & x\geq 0, \\ u_2(t,-x), & x\leq 0, \end{array}\right.,
\end{equation}
and
\begin{equation}
\lambda(x)=\left\{ \begin{array}{ll} \lambda_1(x), & x\geq 0, \\ \lambda_2(-x), & x\leq 0, \end{array}\right.,
\end{equation}
so that
\begin{eqnarray}
u_{1t}&=&\epsilon u_{1xx}+\lambda_{1}(x) u_1\\
u_{2t}&=&\epsilon u_{2xx}+\lambda_{2}(x) u_2,\\
u_1(t,L)&=&U_1(t),\quad
u_2(t,L)=U_2(t),\\
u_1(t,0)&=&u_2(t,0),\quad u_{2x}(t,0)=-u_{1x}(t,0), \label{eqn-bcfold}
\end{eqnarray}
and the transformation as
\begin{eqnarray}
w_1(x)&=&u_1(x)-\int_{0}^{x} K_{11}(x,\xi) u_1(\xi)d\xi
\nonumber \\ &&
-\int_{0}^{x} K_{12}(x,\xi) u_2(\xi)d\xi,\\
w_2(x)&=&u_2(x)-\int_{0}^{x} K_{21}(x,\xi) u_1(\xi)d\xi
\nonumber \\ &&
-\int_{0}^{x} K_{22}(x,\xi) u_2(\xi)d\xi.
\end{eqnarray}
The resulting kernels would be the four ``pieces'' of the kernel $K$ in the transformation (\ref{eqn-transf-unidim}):
\begin{equation}
K(x,\xi)=\left\{\begin{array}{ll}K_{11}(x,\xi), & x\geq 0, 0\leq \xi \leq x ,\\
K_{12}(x,-\xi), & x\geq 0, 0\geq \xi \geq -x , \\
K_{21}(-x,\xi), & x\leq 0, 0\leq \xi \leq -x,\\
K_{22}(-x,-\xi), &x\leq 0, 0\geq \xi \geq x.
 \end{array} \right.
\end{equation}
Therefore the methods of~\cite{orlov}---with some modifications to account for boundary conditions~(\ref{eqn-bcfold})---can potentially be applied to solve reaction-diffusion bilateral control problems.
\section{One-dimensional $2\times2$ hyperbolic linear PDEs with same transport speeds}\label{sect-hyp}
Consider the following system
\begin{eqnarray}
u_t&=&-\epsilon u_x +c_1(x) u+c_2(x) v,\label{eqn-u}\\
v_t&=&\epsilon v_x+c_3(x) u+c_4(x) v ,\label{eqn-v}
\end{eqnarray}
evolving in $x\in[-L,L]$, and assume that $\epsilon>0$ and the coefficients $c_i(x)$ are differentiable. The boundary conditions are:
\begin{eqnarray}
u(t,-L)&=&U_1(t),\quad
v(t,L)=U_2(t),\label{eqn-contr}
\end{eqnarray}
and we need to find $U_1(t)$ and $U_2(t)$ to stabilize the system. 

As in Section~\ref{sect-parabolic} the approach to design $U_1(t)$ and $U_2(t)$ is to seek a mapping that transforms (\ref{eqn-u})--(\ref{eqn-v}) into a target system with adequate properties, in this case
\begin{eqnarray}
\alpha_t&=&-\epsilon \alpha_x+c_1(x)\alpha,\label{eqn-alpha}\\
\beta_t&=&\epsilon \beta_x+c_4(x)\beta ,\label{eqn-beta}
\end{eqnarray}
with boundary conditions
\begin{eqnarray}
\alpha(t,-L)&=&
\beta(t,L)=0,
\end{eqnarray}
whose origin is exponentially stable (in fact it converges to zero in finite time),
and then $U_1(t)$, $U_2(t)$ will be chosen to realize the transformation.

\subsection{Backstepping transformation and kernel equations}
To transform the original system (\ref{eqn-u})--(\ref{eqn-v}) into the target system (\ref{eqn-alpha})--(\ref{eqn-beta}), we look for a mapping defined as follows:
\begin{eqnarray}
\alpha(t,x)&=&u(t,x)-\int_{-x}^x K^{uu}(x,\xi) u(t,\xi) d\xi
\nonumber \\ &&
-\int_{-x}^x K^{uv}(x,\xi) v(t,\xi) d\xi,\label{eqn-tranu}\\
\beta(t,x)&=&v(t,x)-\int_{-x}^x K^{vu}(x,\xi) u(t,\xi) d\xi
\nonumber \\ &&
-\int_{-x}^x K^{vv}(x,\xi) v(t,\xi) d\xi.\label{eqn-tranv}\end{eqnarray}
Introducing (\ref{eqn-tranu})--(\ref{eqn-tranv}) into (\ref{eqn-alpha})--(\ref{eqn-beta}), one obtains the equations that the kernels must satisfy. Defining
\begin{eqnarray}
K(x,\xi)&=&
\left(
\begin{array}{cc}
K^{uu}(x,\xi)  & K^{uv}(x,\xi)     \\
K^{vu}(x,\xi)  & K^{vv}(x,\xi)     
\end{array}
\right),
\Sigma=\epsilon \mathrm{I},\\
C(x)&=&
\left(
\begin{array}{cc}
 c_1(x)  &c_2(x)   \\
c_3(x)  & c_4(x)
\end{array}
\right)\,\\
D(x)&=&
\left(
\begin{array}{cc}
 c_1(x) & 0     \\
0  & c_4(x)
\end{array}
\right)\,\\
 w(t,x)&=&
\left(
\begin{array}{c}
u(t,x) \\
v(t,x)
\end{array}
\right)\,\,,
\gamma (t,x)=
\left(
\begin{array}{c}
\alpha(t,x) \\
\beta(t,x)
\end{array}
\right).\quad
\end{eqnarray}
Then the original plant, target system and transformation can be written compactly as
\begin{eqnarray}
 w_t&=&\Sigma w_x+C(x) w,\\
 w(t,-L) &=&\left( \begin{array}{c}  U_1(t)\\0 \end{array} \right),\,  w(t,L) =\left( \begin{array}{c} 0 \\ U_2(t) \end{array} \right), \\
\gamma_t&=&\Sigma \bold \gamma_x+D(x) \bold \gamma,\\
 \bold \gamma(t,-L)&=& \bold \gamma (t,L) =0\\
\gamma(t,x)&=& w(t,x)-\int_{-x}^x K(x,\xi)  w(t,\xi)d\xi.
\end{eqnarray}
Introducing the transformation into the target system we find a set of three matrix equations:
\begin{eqnarray}
0&=&\Sigma K_x+K_\xi \Sigma-K C(\xi)+D(x)K,\\
0&=&C(x)-D(x)+\Sigma K(x,x)-K(x,x)\Sigma,\\
0&=&\Sigma K(x,-x)-K(x,-x)\Sigma.
\end{eqnarray}
Expanding these terms we get two uncoupled $2\times2$ systems of hyperbolic 1-D equations.

First, for $K^{uu}$ and $K^{uv}$:
\begin{eqnarray} \label{eqn-kuu}
 K_x^{uu}+ K_\xi^{uu}&=&\frac{c_1(\xi)-c_1(x)}{\epsilon}K^{uu}+
\frac{c_3(\xi)}{\epsilon} K^{uv},\quad\\
 K_x^{uv}-K_\xi^{uv}&=&\frac{c_4(\xi)-c_1(x)}{\epsilon}K^{uv}
+\frac{c_2(\xi)}{\epsilon} K^{uu} \label{eqn-kuv},\quad\end{eqnarray}
with boundary conditions
\begin{eqnarray} 
 K^{uu}(x,-x)&=&0,
 \label{eqn-bc1}\\
K^{uv}(x,x)&=&-\frac{c_2(x)}{2\epsilon}. \label{eqn-bc2}
\end{eqnarray}

\begin{figure}[ht]
\includegraphics[width=7cm]{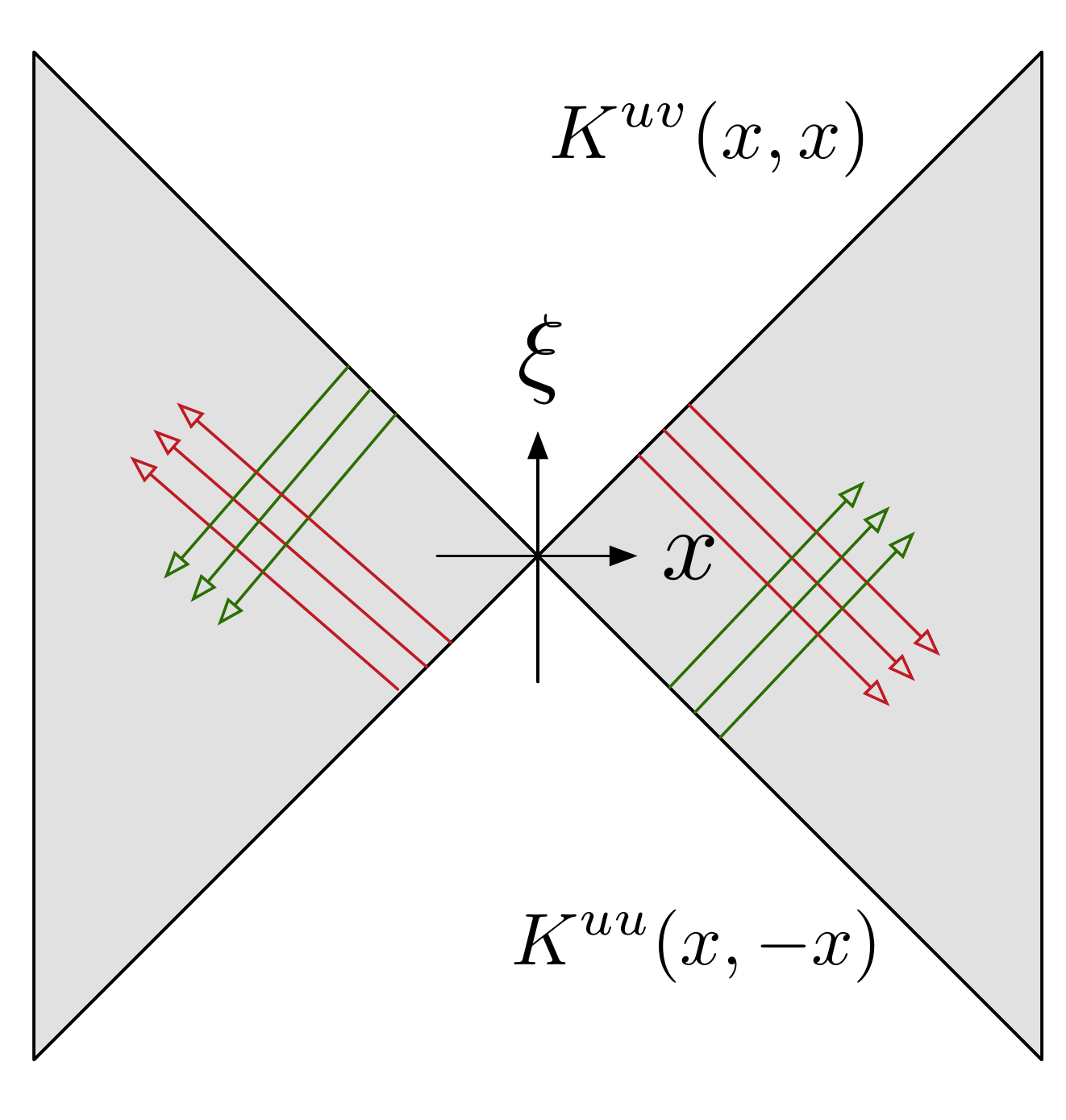}
\centering
\caption{Characteristic lines for the kernel equations. Red: characteristic lines for $K^{uv}$. Green: characteristic lines for $K^{uu}$.}
\label{fig:hourglass2}
\end{figure}

Second, for $K^{vu}$ and $K^{vv}$:
\begin{eqnarray}\label{eqn-kvv}
 K_x^{vv}+K_\xi^{vv}&=&\frac{c_4(x)-c_4(\xi)}{\epsilon}K^{vu}-
\frac{c_2(\xi)}{\epsilon} K^{vu},\quad\\
 K_x^{vu}-K_\xi^{vu}&=&\frac{c_4(x)-c_1(\xi)}{\epsilon}K^{vv}-
\frac{c_3(\xi)}{\epsilon} K^{vv},\quad\end{eqnarray}
with boundary conditions
\begin{eqnarray}
 K^{vv}(x,-x)&=&0,
 \label{eqn-bc3}\\
K^{vu}(x,x)&=&\frac{c_3(x)}{2\epsilon}, \label{eqn-bc4}
\end{eqnarray}
Both systems of equations evolve separately in the domain $\mathcal T$, shown in Figure~\ref{fig:hourglass}. Since they are structurally equivalent, it suffices to analyze one of them, for instance (\ref{eqn-kuu})--(\ref{eqn-bc2}), and as in Section~\ref{sect-parabolic} it is enough to consider just the subdomain $\mathcal T_1$. It is easy to see that the problem has a solution because the boundaries are characteristic (see fig.~\ref{fig:hourglass2}).

Define $y=x+\xi \geq 0$ and $z=x-\xi \geq 0$, and $K(x,\xi)=G(y,z)$. Then the kernel equations become
\begin{eqnarray} \label{eqn-Guu}
 G_y^{uu}&=&
 \frac{c_1\left(\frac{y-z}{2}\right)-c_1\left(\frac{y+z}{2}\right)}{2\epsilon} G^{uu}
\nonumber\\&& +
\frac{c_3\left(\frac{y-z}{2}\right)}{2\epsilon} G^{uv}\\
 G_z^{uv}&=& \frac{c_4\left(\frac{y-z}{2}\right)-c_1\left(\frac{y+z}{2}\right)}{2\epsilon} G^{uv}
\nonumber\\&&+\frac{c_2\left(\frac{y-z}{2}\right)}{2\epsilon} G^{uu} \label{eqn-Guv},\end{eqnarray}
with boundary conditions
\begin{eqnarray} 
G^{uu}(0,z)&=&0
 \label{eqn-bcG1}\\
G^{uv}(y,0)&=&-\frac{c_2(2y)}{2\epsilon}. \label{eqn-bcG2}
\end{eqnarray}
The equations are easily converted into integral equations, namely
\begin{eqnarray} \label{eqn-kuu_integrated22}
 G^{uu}&=&
  \int_0^{y} \frac{c_1\left(\frac{s-z}{2}\right)-c_1\left(\frac{s+z}{2}\right)}{2\epsilon} G^{uu}(s,z) ds
 \nonumber \\ &&
 +\int_0^{y}
\frac{c_3\left(\frac{s-z}{2}\right)}{2\epsilon} G^{uv}(s,z) ds,\\
G^{uv}&=&
\int_0^z \frac{c_4\left(\frac{y-s}{2}\right)-c_1\left(\frac{y+s}{2}\right)}{2\epsilon} G^{uv}(y,s) ds
 \nonumber \\ &&-
\frac{c_2(2y)}{2\epsilon}
+
\int_0^z
\frac{c_2\left(\frac{y-s}{2}\right)}{2\epsilon} G^{uu}(y,s)ds.\quad
\end{eqnarray}
Using a successive approximation series to compute $G^{uu}$ and $G^{uv}$, we get
\begin{equation}
G^{uv}(y,z)=\sum_0^\infty F^{uv}_i(y,z),\quad G^{uu}(y,z)=\sum_0^\infty F^{uu}_i(y,z),
\end{equation}
with $F^{uu}_0=0$,
\begin{eqnarray} 
F_0^{uv}
&=&-\frac{c_2(2y)}{2\epsilon},\\
F_i^{uu}&=&
\int_0^{y} \frac{c_1\left(\frac{s-z}{2}\right)-c_1\left(\frac{s+z}{2}\right)}{2\epsilon} F_{i-1}^{uu}(s,z) ds
 \nonumber \\ &&
 +\int_0^{y}
\frac{c_3\left(\frac{s-z}{2}\right)}{2\epsilon} F_{i-1}^{uv}(s,z) ds, \\
F_i^{uv}&=&
\int_0^z \frac{c_4\left(\frac{y-s}{2}\right)-c_1\left(\frac{y+s}{2}\right)}{2\epsilon}  F_{i-1}^{uv}(y,s) ds
 \nonumber \\ &&
+
\int_0^z
\frac{c_2\left(\frac{y-s}{2}\right)}{2\epsilon}  F_{i-1}^{uu}(y,s)ds.
\end{eqnarray}
Calling $\lambda=\dfrac{1}{2\epsilon}\max_{x\in[-L,L]}\left\{c_1(x),c_2(x),c_3(x),c_4(x)\right\}$, let us show that the very simple bounds
\begin{equation}
\vert F_i^{uu}(y,z) \vert , \vert F_i^{uv}(y,z) \vert \leq 4^i \lambda^{i+1} \frac{(y+z)^i}{i!}
\end{equation}
works. Obviously it does for $i=0$. Now assuming it is correct for $i-1$, we obtain, for instance for $F_i^{uu}$,
\begin{eqnarray} 
\vert F_i^{uu}(y,z) \vert &\leq& 2\lambda 
\int_0^{y} \left(\vert F_{i-1}^{uu}(s,z) \vert+\vert F_{i-1}^{uv}(s,z)\vert\right) ds \nonumber  \\
&\leq &
\frac{ 4^{i} \lambda^{i+1}}{(i-1)!}
\int_0^{y} (s+z) ds \nonumber  \\
&\leq & 4^i \lambda^{i+1} \frac{(y+z)^i}{i!},
\end{eqnarray}
and similarly for $F_i^{uv}$, thus proving the bound.
Therefore,
\begin{eqnarray}
\vert G^{uv}(y,z)\vert 
&\leq& \sum_0^\infty\vert F^{uv}_i(y,z) \vert
\nonumber \\Ê
&\leq& \sum_0^\infty 4^i \lambda^{i+1} \frac{(y+z)^i}{i!}
\nonumber \\Ê
&=& \lambda \mathrm{e}^{4\lambda(y+z)},
\end{eqnarray}
and similarly for $G^{uu}$. Therefore, extending the proof to $\mathcal T_2$,
\begin{equation}
\vert K^{uv}(x,\xi)\vert \leq  \lambda \mathrm{e}^{8\lambda \vert x\vert},
\end{equation}
and the same bound for $K^{uu}$, $K^{vv}$, and $K^{vu}$ applies.

The resulting feedback laws are:
\begin{eqnarray}
U_1(t)&=&\int_{-L}^L K^{uu}(-L,\xi) u(t,\xi) d\xi \nonumber \\Ê&&
+\int_{-L}^L K^{uv}(-L,\xi) v(t,\xi) d\xi,\label{eqn-conu}\\
U_2(t)&=&\int_{-L}^L K^{vu}(L,\xi) u(t,\xi) d\xi\nonumber \\Ê&&
+\int_{-L}^L K^{vv}(L,\xi) v(t,\xi) d\xi,\label{eqn-conv}\end{eqnarray}
\section{Wave equation}\label{sect-wave}
Consider the following hyperbolic PDE
\begin{equation}
u_{tt}-u_{xx}=2 \lambda(x) u_t+\alpha(x) u_x+\beta(x)u,
\end{equation}
in the domain $x\in[-L,L]$, which represents a wave equation with (potentially) in-domain anti-damping~\cite{krstic2},
with boundary conditions
\begin{equation}
u(t,-L)=U_1(t),\qquad u(t,L)=U_2(t).
\end{equation}

The idea pursued in this paper to stabilize the wave equation is to identify this equation with a hyperbolic $2\times2$ system. We consider the case\footnote{If $\beta(x)\neq0$ then, depending on the values of $\alpha,\beta,\lambda$ a scaling transformation~\cite{krstic2} may exist that transforms the system into another with $\beta(x)=0$. If not, then using an alternate definition for the states $w$ and $v$ the problem can still be solved. We leave the details outfor lack of space.} in which $\beta(x)=0$.

Define $w = u_x + u_t$ and $v = u_x - u_t$. Notice that $v+w=2u_x$ and $w-v=2u_t$. Then 
\begin{eqnarray}\label{eqn-wwave}
w_t-w_x&=&
\lambda(x)(w-v)+\frac{\alpha(x)}{2} (w+v),\\
v_t+v_x&=&
\lambda(x)(v-w)-\frac{\alpha(x)}{2} (w+v).
\end{eqnarray}
To find the boundary conditions, we notice
\begin{eqnarray}
w(-1)&=&u_x(t,-1)+\dot U_1=V_1,\\
v(1)&=&u_x(t,1)-\dot U_2=V_2,\label{eqn-bc2wave}
\end{eqnarray}
where $V_1$ and $V_2$ have been artificially introduced. Since (\ref{eqn-wwave})--(\ref{eqn-bc2wave}), with the coefficients $c_i(x)$ found as combinations of $\alpha$ and $\lambda$, one can find the values of $V_1$ and $V_2$ by applying the design of Section~\ref{sect-hyp}. Then, the values of $U_1$ and $U_2$ are found from solving the differential equation
\begin{eqnarray}
\dot U_1&=&V_1-u_x(t,-1)\\
\dot U_2&=&-V_2+u_x(t,1)
\end{eqnarray}
with initial conditions for $U_1$ and $U_2$ consistent with, respectively, $u(t,-1)$ and $u(t,1)$.
\section{Explicit control laws}\label{sect-explicit}
In this section we provide explicit formulae for the control problems posed in this paper, in the constant coefficient case. These expressions are subsequently used to compare unilateral and bilateral control laws.
\subsection{Reaction-Diffusion Equation}
In the case of constant $\lambda$, the kernel equations (\ref{eqn-K1})--(\ref{eqn-K3}) are
\begin{eqnarray} \label{Kexact1}
 K_{xx}(x,\xi)- K_{\xi \xi}(x,\xi)&=&
\frac{\lambda}{\epsilon} K(x,\xi),\\
K(x,x)&=&- \frac{\lambda}{2\epsilon}x ,\\
K(x,-x)&=&0.\label{Kexact3}
\end{eqnarray}
Using the techniques of~\cite{nball} we find a explicit solution as
\begin{equation}\label{eqn-kernel-unidim}
K(x,\xi)=-\mathrm{sgn}(x)\frac{1}{2}\sqrt{\frac{\lambda}{\epsilon}} \mathrm{I}_1 \left[\sqrt{\frac{\lambda}{\epsilon} \left(x^2-\xi^2\right)}\right] \sqrt{\frac{x+\xi}{x-\xi}}.
\end{equation}
Thus, the following control laws stabilize the system:
\begin{eqnarray}\label{eqn-1dim-controlaw-rect}
U_1\hspace{-7pt}
&=&\hspace{-7pt}-\frac{1}{2}\sqrt{\frac{\lambda}{\epsilon}} \int_{-L}^L \hspace{-2pt}\sqrt{\frac{L+\xi}{L-\xi}} \mathrm{I}_1\hspace{-3pt}\left[\sqrt{\frac{\lambda}{\epsilon}(L^2-\xi^2)}\right] \hspace{-3pt}  u(t,\xi) d \xi,\quad\,\,\,\,\\
U_2
\hspace{-7pt}
&=&\hspace{-7pt}- \frac{1}{2}\sqrt{\frac{\lambda}{\epsilon}} \int_{-L}^L  \hspace{-2pt}\sqrt{\frac{L-\xi}{L+\xi}} \mathrm{I}_1\hspace{-3pt}\left[\sqrt{\frac{\lambda}{\epsilon}(L^2-\xi^2)}\right]  \hspace{-3pt} u(t,\xi) d \xi.\quad\,\,\,\,\label{eqn-1dim-controlaw-rect2}
\end{eqnarray}
\subsection{One-dimensional $2\times2$ hyperbolic linear PDEs with same transport speeds}
In the case of constant coefficients, and scaling the kernels by the function $\exp\left(\dfrac{(c_4-c_1)(x-\xi)}{2\epsilon}\right) $, it is possible to reduce (\ref{eqn-kuu})--(\ref{eqn-bc4}) to (\ref{Kexact1})--(\ref{Kexact3}), with $\lambda=c_2c_3$. We leave the details out for lack of space. The reached solution is
\begin{eqnarray}
 K^{uu}&=&\exp\left[\dfrac{(c_4-c_1)(x-\xi)}{2\epsilon}\right]  F(x,\xi),\\
 K^{vv}&=&\exp\left[\dfrac{(c_4-c_1)(x-\xi)}{2\epsilon}\right]  F(x,\xi),\\
K^{uv}&=& \frac{c_2}{2\epsilon}\exp\left[\dfrac{(c_4-c_1)(x-\xi)}{2\epsilon}\right]  H(x,\xi),\\
K^{vu}&=&- \frac{c_3}{2\epsilon}\exp\left[\dfrac{(c_4-c_1)(x-\xi)}{2\epsilon}\right]  H(x,\xi),
\end{eqnarray}
where
\begin{eqnarray}
 F\hspace{-3pt}&\hspace{-3pt}=\hspace{-3pt}&\hspace{-3pt}-\mathrm{sgn}(x)\frac{\sqrt{c_2c_3}}{2\epsilon}\sqrt{\frac{x+\xi}{x-\xi}} \mathrm{I}_1\left[\frac{ \sqrt{c_2c_3(x^2-\xi^2)}}{\epsilon}\right],\quad\\
 H \hspace{-3pt}&\hspace{-3pt}=\hspace{-3pt}&\hspace{-3pt}-\mathrm{sgn}(x)  \mathrm{I}_0\left[\frac{ \sqrt{c_2c_3(x^2-\xi^2)}}{\epsilon}\right].
\end{eqnarray}
The control laws are then
\begin{eqnarray}
U_1&=&\int_{-L}^L
\left\{\frac{c_2}{2\epsilon} \mathrm{I}_0\left[\frac{ \sqrt{c_2c_3(L^2-\xi^2)}}{\epsilon}\right] v(t,\xi)
  \right. \nonumber \\ &&
\left.+ \frac{\sqrt{c_2c_3}}{2\epsilon}\sqrt{\frac{L-\xi}{L+\xi}} \mathrm{I}_1\left[\frac{ \sqrt{c_2c_3(L^2-\xi^2)}}{\epsilon}\right] u(t,\xi)
\right\}
\nonumber \\Ê&& \times
\exp\left[\dfrac{(c_4-c_1)(L-\xi)}{2\epsilon}\right] 
 d\xi,\label{eqn-conu2}\\
U_2&=&
\int_{-L}^L
\left\{\frac{c_3}{2\epsilon} \mathrm{I}_0\left[\frac{ \sqrt{c_2c_3(L^2-\xi^2)}}{\epsilon}\right] u(t,\xi) 
\right.
\nonumber \\ &&
-\left.  \frac{\sqrt{c_2c_3}}{2\epsilon}\sqrt{\frac{L+\xi}{L-\xi}} \mathrm{I}_1\left[\frac{ \sqrt{c_2c_3(L^2-\xi^2)}}{\epsilon}\right]  v(t,\xi) \right\}\nonumber \\Ê&& \times
\exp\left[\dfrac{(c_4-c_1)(L-\xi)}{2\epsilon}\right] 
 d\xi
.\label{eqn-conv2}\end{eqnarray}

\section{Comparison between unilateral and bilateral feedback laws}\label{sect-comparison}
Unilateral and bilateral feedback laws can be compared  by setting one of the controllers to zero (for instance $U_1=0$) and using the single-sided control design. We concentrate on constant-coefficient cases because the explicit laws are available and simplify the comparison and only consider reaction-diffusion equations. The results are similar for the other equations and they will be more thoroughly studied in a future journal publication.

\begin{figure}[ht]
\includegraphics[width=7.5cm]{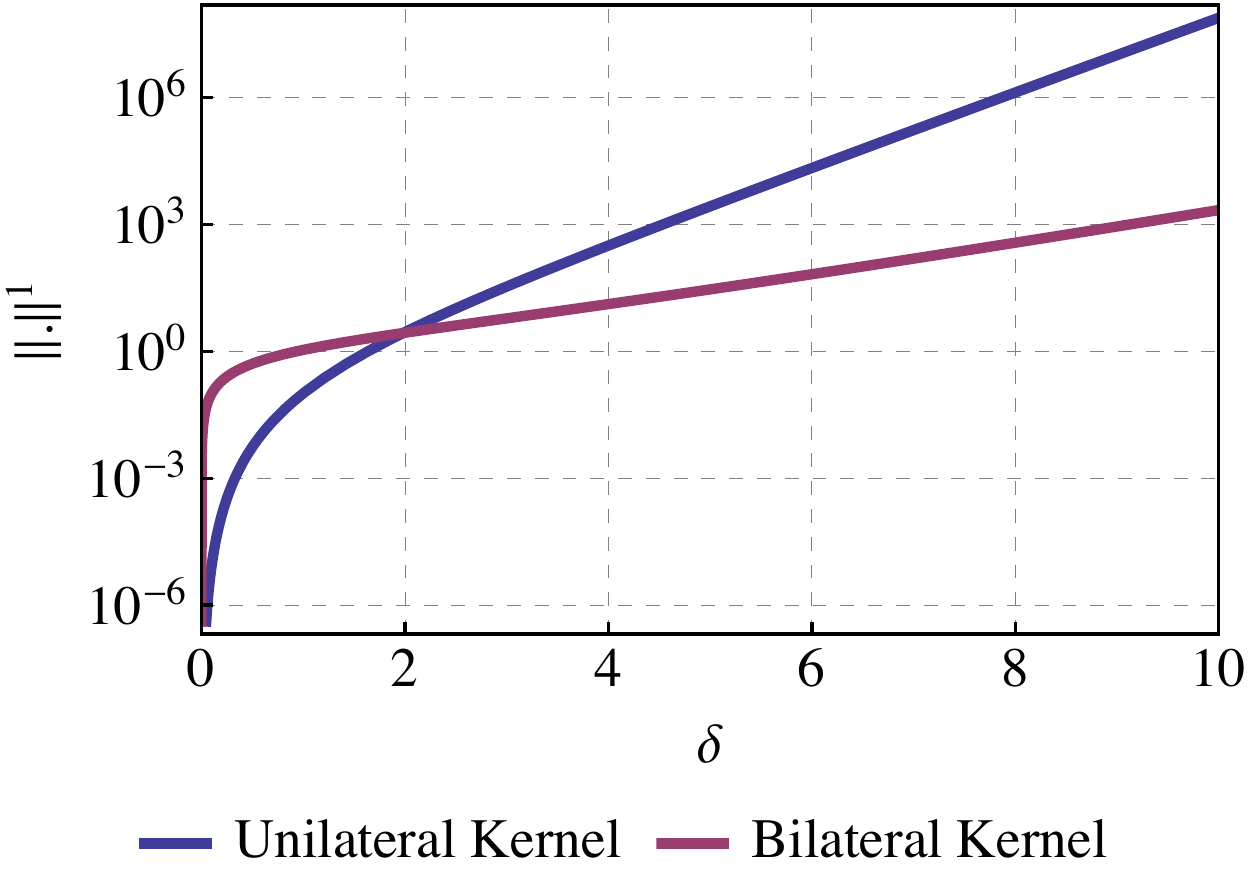}
\centering
\caption{Plot of unilateral and bilateral $L^1$ norm of reaction-diffusion control kernels for varying $\delta=L\sqrt{\frac{\lambda}{\epsilon}}$.}
\label{fig:compkernelsRD}
\end{figure}

For reaction-diffusion equations, the control law is~\cite{krstic}
\begin{equation} \label{eqn-controluni}
U_2=\int_{-L}^L -\sqrt{\frac{\lambda}{\epsilon}} \xi \frac{\mathrm{I}_1 \left[\sqrt{\dfrac{\lambda}{\epsilon}(4L^2-(\xi+L)^2)}\right] }{\sqrt{4L^2-(\xi+L)^2}} u(t,\xi) d\xi.
\end{equation}
where it has been taken into account that the lenght of the domain is $2L$. To compare (\ref{eqn-controluni}) and (\ref{eqn-1dim-controlaw-rect})--(\ref{eqn-1dim-controlaw-rect2}), one can use the $L^1$ norm of the kernels. Call $J_1$ and $J_2$ the $L^1$ norm of the unilateral and bilateral kernels, respectively. Notice that the $L^1$ norm of the kernel in (\ref{eqn-1dim-controlaw-rect}), changing $\xi$ to $L-\xi$, is the same as the kernel in (\ref{eqn-1dim-controlaw-rect2}). Then, one has
\begin{eqnarray}
J_1&=&\sqrt{\frac{\lambda}{\epsilon}} \int_{-L}^L \xi  \frac{\mathrm{I}_1 \left[\sqrt{\dfrac{\lambda}{\epsilon}(4L^2-(\xi+L)^2)}\right] }{\sqrt{4L^2-(\xi+L)^2}} d\xi,\\
J_2&=&\sqrt{\frac{\lambda}{\epsilon}} \int_{-L}^L \sqrt{\frac{L+\xi}{L-\xi}} \mathrm{I}_1\left[\sqrt{\frac{\lambda}{\epsilon}(L^2-\xi^2)}\right] d\xi,
\end{eqnarray}
which,  writing $\delta=L\sqrt{\frac{\lambda}{\epsilon}}$, can be expressed as
\begin{eqnarray}
J_1&=&\delta \int_{-1}^1  \xi\frac{\mathrm{I}_1 \left[\delta\sqrt{4-(\xi+1)^2}\right] }{\sqrt{4-(\xi+1)^2}} d\xi,\\
J_2&=&\delta \int_{-1}^1 \sqrt{\frac{1+\xi}{1-\xi}} \mathrm{I}_1\left[\delta\sqrt{1-\xi^2}\right] d\xi,\quad
\end{eqnarray}
which is exclusively a function of $\delta$. In Fig.~\ref{fig:compkernelsRD} $J_1$ and $J_2$ are represented for different values of $\delta$. It can be seen that for smaller values of $\delta$ (in concrete, $\delta<2$), the unilateral kernel is smaller, however for larger values of $\delta$ ($\delta>2$), the bilateral kernel grows more slowly.

Besides a comparison in terms of magnitude, it is clear that a bilateral design can be made fault-tolerant, in the sense that it can withstand the loss of one of the actuators. Once the failure is detected, it is sufficient to switch to a unilateral control law with the remaining controller. Detecting the failure would require the use of an additional observer, which would require at least one measurement sensor. Fault-tolerant designs will be explored in future works.
\section{Concluding remarks}\label{sect:conclusions}

This paper presents bilateral control laws for linear reaction-diffusion, wave and $2\times2$ hyperbolic 1-D systems. For $2\times2$ hyperbolic 1-D systems, only the case of both states having the same speed of transport is considered. The control design for the wave equation is built upon the design for hyperbolic systems. The unilateral and bilateral control laws are compared for reaction-diffusion equations, showing that the bilateral control law requires less total effort (despite using two actuators) when the system coefficients are large.

These results open the door for more sophisticated designs such as bilateral sensor/actuator output feedback and fault-tolerant designs that will be explored in future publications.

\section*{Acknowledgments}

Rafael Vazquez acknowledges financial
support of the Spanish Ministerio de Econom\'ia y Competitividad  under grant MTM2015-65608-P.

\end{document}